\documentclass[review]{elsarticle}

\usepackage{bm,amssymb,amsmath,mathrsfs}
\usepackage{amsthm}
\usepackage{lineno}
\usepackage[usenames]{color}

\usepackage{dcolumn}

\newcommand{\bq}{\begin{equation}}
\newcommand{\eq}{\end{equation}}
\newcommand{\bc}{\begin{center}}
\newcommand{\ec}{\end{center}}
\newcommand{\bit}{\begin{itemize}}
\newcommand{\eit}{\end{itemize}}
\newcommand{\ben}{\begin{enumerate}}
\newcommand{\een}{\end{enumerate}}

\theoremstyle{plain}

\newtheorem*{theorem*}{Theorem}

\begin{document}

\journal{(internal report CC23-9)}

\begin{frontmatter}

\title{Alternative combinatorial sum for the probability mass function of the Poisson distribution of order $k$}

\author[cc]{S.~R.~Mane}
\ead{srmane001@gmail.com}
\address[cc]{Convergent Computing Inc., P.~O.~Box 561, Shoreham, NY 11786, USA}

\begin{abstract}
  Kostadinova and Minkova published an expression for the probability mass function (pmf) of the Poisson distribution of order $k$,
  as a combinatorial sum ($\mathit{Pliska~Stud.~Math.~Bulgar.}\ {\bf 22},\ 117-128\ (2013)$).
  Inspired by their elegant solution, this note presents an alternative combinatorial sum for the pmf of the Poisson distribution of order $k$.
  The terms are partitioned into blocks of length $k$ (as opposed to $k+1$ by Kostadinova and Minkova).
  The new sum offers an advantage in the following sense.
  For $n\in[rk+1,(r+1)k]$, the lowest power of $\lambda$ in the pmf is $\lambda^{r+1}$.
  Hence the lower limit of summation can be increased, to avoid needlessly calculating terms which cancel to identically zero.
\end{abstract}

\vskip 0.25in


\end{frontmatter}

\newpage
The Poisson distribution of order $k$ is a special case of a compound Poisson distribution introduced by Adelson \cite{Adelson1966}.
For a (possibly infinite) tuple $\bm{a}=(a_1,a_2,\dots)$ and $x\in\mathbb{R}$, the probability generating function (pgf)
of a compound Poisson distribution is (eq.~(1) in \cite{Adelson1966})
\bq
\label{eq:Adelson_pgf}
f(\bm{a},x) = \exp\Bigl(-\sum_i a_i\Bigr)\exp\Bigl(\sum_i a_ix^i\Bigr) \,.
\eq
The Poisson distribution of order $k$ is the special case where $a_1=\dots=a_k=\lambda$, where $\lambda>0$, and all the other $a_i$ are zero.
For $k=1$ it is the standard Poisson distribution.
The probability mass function (pmf) of the Poisson process of order $k$ is given as follows (eq.~(1) in \cite{KwonPhilippou}, with slight changes of notation)
\bq
\label{eq:KP_pmf}
f_k(n;\lambda) = e^{-k\lambda} \sum_{n_1+2n_2+\dots+kn_k=n} \frac{\lambda^{n_1+\dots+n_k}}{n_1!\dots n_k!} \qquad (n=0,1,2,\dots)\,.
\eq
Observe that the sum is a polynomial in $\lambda$.
The following facts are easily deduced from eq.~\eqref{eq:KP_pmf}:
\begin{enumerate}
\item
  For $n=0$, the sum equals $1$.
\item
  For fixed $n>0$, the sum does not have a constant term.
\item
  For fixed $n>0$, the highest power of $\lambda$ in the sum is $\lambda^n$, i.e.~the sum is a polynomial of degree $n$.
\item
  If $n\in[rk+1,(r+1)k]$, where $r=0,1,\dots$, the {\em lowest} power of $\lambda$ in the sum is $\lambda^{r+1}$.
\end{enumerate}
For most of this note, we shall hold $k\ge2$ and $\lambda>0$ fixed and vary only the value of $n$.
We adopt the notation by Kostadinova and Minkova \cite{KostadinovaMinkova2013}
and write ``$p_n$'' in place of $f_k(n;\lambda)$ and omit explicit mention of $k$ and $\lambda$. 
Kostadinova and Minkova published the following expression for the pmf, in terms of combinatorial sums
(Theorem 1 in \cite{KostadinovaMinkova2013}, with slight changes of notation)
\bq
\label{eq:KM_Thm1}
\begin{split}
  p_0 &= e^{-k\lambda} \,,
  \\
  p_n &= e^{-k\lambda}\,\sum_{j=1}^n \binom{n-1}{j-1}\,\frac{\lambda^j}{j!} \qquad (n=1,2,\dots,k) \,,
  \\
  p_n &= e^{-k\lambda}\,\biggl[\,\sum_{j=1}^n \binom{n-1}{j-1}\,\frac{\lambda^j}{j!}
  \;\;-\;\; \sum_{i=1}^\ell (-1)^{i-1}\,\frac{\lambda^i}{i!} \sum_{j=0}^{n-i(k+1)} \binom{n - i(k+1) +i-1}{j+i-1}\,\frac{\lambda^j}{j!} \,\biggr]
  \\
  &\qquad\qquad (n = \ell(k+1)+m,\, m=0,1,\dots,k,\, \ell=1,2,\dots,\infty) \,.
\end{split}
\eq
The sum in eq.~\eqref{eq:KM_Thm1}, although elegant, contains terms which cancel to identically zero, for low powers of $\lambda$.
Observe also that the tail piece for $n>k$ is composed of blocks of $k+1$ elements each, for $\ell=1,2,\dots$.
This makes it difficult to identify the terms in eq.~\eqref{eq:KM_Thm1} which cancel to identically zero.
The author was therefore motivated to find an alternative combinatorial sum, where the terms for $n>k$ are partitioned into blocks of $k$ elements
$n\in[rk+1,(r+1)k]$, where $r=0,1,\dots$.
Then $r = \lfloor(n-1)/k\rfloor$.
The result is
\bq
\label{eq:KM_pmf_me}
\begin{split}
  p_0 &= e^{-k\lambda} \,,
  \\
  p_n &= e^{-k\lambda}\,\sum_{j=1}^n \binom{n-1}{j-1}\,\frac{\lambda^j}{j!} \qquad (n=1,2,\dots,k) \,,
  \\
  p_n &= e^{-k\lambda}\,\biggl[\, \sum_{j=r+1}^n \binom{n-1}{j-1}\,\frac{\lambda^j}{j!}
  \;-\; \sum_{i=1}^r (-1)^{i-1}\,\frac{\lambda^i}{i!} \sum_{j=r+1-i}^{n-ik-1} \binom{n - ik -1}{j+i-1}\,\frac{\lambda^j}{j!} \,\biggr]
  \\
  &\qquad\qquad (n = rk+m,\, m=1,\dots,k,\, r=1,2,\dots,\infty) \,.
\end{split}
\eq
There is no change relative to eq.~\eqref{eq:KM_Thm1} for $n \le k$.
Observe that for $n>k$, the lower limit of the sums do {\em not} begin from $j=1$ (as they do in eq.~\eqref{eq:KM_Thm1}).
Terms which cancel to zero do not appear in eq.~\eqref{eq:KM_pmf_me}.
Note that in eq.~\eqref{eq:KM_pmf_me}, it is possible (if $n>k$) for the lower limit of the sum to be higher than the upper limit: in such a case the value of the sum is zero.
Let us validate eq.~\eqref{eq:KM_pmf_me} for $k=2$ and a few values of $n$.
(The same exercise was carried out in \cite{Mane_Poisson_k_CC23_8}, for eq.~\eqref{eq:KM_Thm1}, where the cancellation of low powers of $\lambda$ was observed,
which is the motivation for this note.)
As in \cite{Mane_Poisson_k_CC23_8}, for brevity we omit the prefactor of $e^{-k\lambda}$ and calculate only the polynomial in $\lambda$ (without change of notation for $p_n$).
For $k=2$ the expression for $p_n$ is simple, for all $n\ge1$ (see eq.~\eqref{eq:KP_pmf}).
\bq
\label{eq:pmf_k2}
\begin{split}
p_n &= \sum_{j=0}^{\lfloor(n/2)\rfloor} \frac{\lambda^{n-j}}{(n-2j)!j!}
\\
&= \frac{\lambda^n}{n!} +\frac{\lambda^{n-1}}{(n-2)!1!} +\frac{\lambda^{n-2}}{(n-4)!2!} +\dots +\frac{\lambda^{n-\lfloor(n/2)\rfloor}}{\lfloor(n/2)\rfloor!} \,.
\end{split}
\eq
The first few polynomials are as follows.
\begin{subequations}
\label{eq:poly_k2}
\begin{align}
  p_0 &= 1 \,,
  \\
  p_1 &= \lambda \,,
  \\
  p_2 &= \frac{\lambda^2}{2!} + \frac{\lambda}{0!1!} \,,
  \\
  p_3 &= \frac{\lambda^3}{3!} + \frac{\lambda^2}{1!1!} \,,
  \\
  p_4 &= \frac{\lambda^4}{4!} + \frac{\lambda^3}{2!1!} + \frac{\lambda^2}{0!2!} \,,
  \\
  p_5 &= \frac{\lambda^5}{5!} + \frac{\lambda^4}{3!1!} + \frac{\lambda^3}{1!2!} \,,
  \\
  p_6 &= \frac{\lambda^6}{6!} + \frac{\lambda^5}{4!1!} + \frac{\lambda^4}{2!2!} + \frac{\lambda^3}{0!3!} \,,
  \\
  p_7 &= \frac{\lambda^7}{7!} + \frac{\lambda^6}{5!1!} + \frac{\lambda^5}{3!2!} + \frac{\lambda^4}{1!3!} \,,
  \\
  p_8 &= \frac{\lambda^8}{8!} + \frac{\lambda^7}{6!1!} + \frac{\lambda^6}{4!2!} + \frac{\lambda^5}{2!3!} + \frac{\lambda^4}{0!4!} \,.
\end{align}
\end{subequations}
\newpage
\noindent
Let us employ eq.~\eqref{eq:KM_pmf_me} and compare with the expressions in eq.~\eqref{eq:poly_k2}.
\begin{enumerate}
\item
  Case $n=1$.
\bq
\begin{split}
  p_1 &= \sum_{j=1}^1 \binom{1-1}{j-1}\,\frac{\lambda^j}{j!}
  = \lambda \,.
\end{split}
\eq
\item
  Case $n=2$.
\bq
\begin{split}
  p_2 &= \sum_{j=1}^2 \binom{2-1}{j-1}\,\frac{\lambda^j}{j!}
  = \binom{1}{0}\lambda +\binom{1}{1}\frac{\lambda^2}{2!}
  = \frac{\lambda}{0!1!} +\frac{\lambda^2}{2!0!} \,.
\end{split}
\eq
\item
  Case $n=3$, then $r=1$.
\bq
\begin{split}
  p_3 &= \sum_{j=2}^3 \binom{2}{j-1}\,\frac{\lambda^j}{j!}
  \;-\; \sum_{i=1}^1 (-1)^{i-1}\,\frac{\lambda^i}{i!} \sum_{j=2-i}^{2-2i} \binom{3 - 2i -1}{j+i-1}\,\frac{\lambda^j}{j!}
  \\
  &= \binom{2}{1}\frac{\lambda^2}{2!} +\binom{2}{2}\frac{\lambda^3}{3!}  - 0
  \\
  &= \frac{\lambda^2}{1!1!} +\frac{\lambda^3}{3!0!} \,.
\end{split}
\eq
\item
  Case $n=4$, then $r=1$.
\bq
\begin{split}
  p_4 &= \sum_{j=2}^4 \binom{3}{j-1}\,\frac{\lambda^j}{j!}
  \;-\; \sum_{i=1}^1 (-1)^{i-1}\,\frac{\lambda^i}{i!} \sum_{j=2-i}^{3-2i} \binom{4 - 2i -1}{j+i-1}\,\frac{\lambda^j}{j!}
  \\
  &= \binom{3}{1}\frac{\lambda^2}{2!} +\binom{3}{2}\frac{\lambda^3}{3!} +\binom{3}{3}\frac{\lambda^4}{4!}
  -\lambda \sum_{j=1}^{1} \binom{1}{j}\,\frac{\lambda^j}{j!}
  \\
  &= \frac{3\lambda^2}{0!2!} +\frac{\lambda^3}{2!1!} +\frac{\lambda^4}{4!0!} -\lambda^2
  \\
  &= \frac{\lambda^2}{0!2!} +\frac{\lambda^3}{2!1!} +\frac{\lambda^4}{4!0!} \,.
\end{split}
\eq

\item
  Case $n=5$, then $r=2$.
\bq
\begin{split}
  p_5 &= \sum_{j=3}^5 \binom{4}{j-1}\,\frac{\lambda^j}{j!}
  \;-\; \sum_{i=1}^2 (-1)^{i-1}\,\frac{\lambda^i}{i!} \sum_{j=3-i}^{4-2i} \binom{5 - 2i -1}{j+i-1}\,\frac{\lambda^j}{j!}
  \\
  &= \binom{4}{2}\frac{\lambda^3}{3!} +\binom{4}{3}\frac{\lambda^4}{4!} +\binom{4}{4}\frac{\lambda^5}{5!}
  -\lambda \sum_{j=2}^{2} \binom{2}{j}\,\frac{\lambda^j}{j!}
  +\frac{\lambda^2}{2!} \sum_{j=1}^{0} \binom{0}{j}\,\frac{\lambda^j}{j!}
  \\
  &= \frac{6\lambda^3}{3!} +\frac{4\lambda^4}{4!} +\frac{\lambda^5}{5!}
  -\frac{\lambda^3}{2!}
  \\
  &= \frac{\lambda^3}{1!2!} +\frac{\lambda^4}{3!1!} +\frac{\lambda^5}{5!0!} \,.
\end{split}
\eq

\item
  Case $n=6$, then $r=2$.
\bq
\begin{split}
  p_6 &= \sum_{j=3}^6 \binom{5}{j-1}\,\frac{\lambda^j}{j!}
  \;-\; \sum_{i=1}^2 (-1)^{i-1}\,\frac{\lambda^i}{i!} \sum_{j=3-i}^{5-2i} \binom{6 - 2i -1}{j+i-1}\,\frac{\lambda^j}{j!}
  \\
  &= \binom{5}{2}\frac{\lambda^3}{3!} +\binom{5}{3}\frac{\lambda^4}{4!}
  +\binom{5}{4}\frac{\lambda^5}{5!} +\binom{5}{5}\frac{\lambda^6}{6!}
  \\
  &\quad
  -\lambda \sum_{j=2}^{3} \binom{3}{j}\,\frac{\lambda^j}{j!}
  +\frac{\lambda^2}{2!} \sum_{j=1}^{1} \binom{1}{j+1}\,\frac{\lambda^j}{j!}
  \\
  &= \frac{10\lambda^3}{3!} +\frac{10\lambda^4}{4!} +\frac{5\lambda^5}{5!} +\frac{\lambda^6}{6!}
  \\
  &\quad
  -\lambda \Bigl(\frac{3\lambda^2}{2!} +\frac{\lambda^3}{3!}\Bigr)
  + 0
  \\
  &= \frac{\lambda^3}{0!3!} +\frac{\lambda^4}{2!2!} +\frac{\lambda^5}{4!1!} +\frac{\lambda^6}{6!0!} \,.
\end{split}
\eq

\item
  Case $n=7$, then $r=3$.
\bq
\begin{split}
  p_7 &= \sum_{j=4}^7 \binom{6}{j-1}\,\frac{\lambda^j}{j!}
  \;-\; \sum_{i=1}^3 (-1)^{i-1}\,\frac{\lambda^i}{i!} \sum_{j=4-i}^{6-2i} \binom{7 - 2i -1}{j+i-1}\,\frac{\lambda^j}{j!}
  \\
  &= \binom{6}{3}\frac{\lambda^4}{4!} +\binom{6}{4}\frac{\lambda^5}{5!} +\binom{6}{5}\frac{\lambda^6}{6!} +\binom{6}{6}\frac{\lambda^7}{7!}
  \\
  &\quad
  -\lambda \sum_{j=3}^{4} \binom{4}{j}\,\frac{\lambda^j}{j!}
  +\frac{\lambda^2}{2!} \sum_{j=2}^{2} \binom{2}{j+1}\,\frac{\lambda^j}{j!}
  -\frac{\lambda^3}{3!} \sum_{j=1}^{2} \binom{0}{j+2}\,\frac{\lambda^j}{j!}
  \\
  &= \frac{20\lambda^4}{4!} +\frac{15\lambda^5}{5!} +\frac{6\lambda^6}{6!} +\frac{\lambda^7}{7!}
  \\
  &\quad
  -\lambda \Bigl(\frac{4\lambda^3}{3!} +\frac{\lambda^4}{4!}\Bigr)
  +0 - 0
  \\
  &= \frac{\lambda^4}{1!3!} +\frac{\lambda^5}{3!2!} +\frac{\lambda^6}{5!1!} +\frac{\lambda^7}{7!0!} \,.
\end{split}
\eq

\item
  Case $n=8$, then $r=3$.
\bq
\begin{split}
  p_8 &= \sum_{j=4}^8 \binom{7}{j-1}\,\frac{\lambda^j}{j!}
  \;-\; \sum_{i=1}^3 (-1)^{i-1}\,\frac{\lambda^i}{i!} \sum_{j=4-i}^{7-2i} \binom{8 - 2i -1}{j+i-1}\,\frac{\lambda^j}{j!}
  \\
  &= \binom{7}{3}\frac{\lambda^4}{4!} +\binom{7}{4}\frac{\lambda^5}{5!} +\binom{7}{5}\frac{\lambda^6}{6!} +\binom{7}{6}\frac{\lambda^7}{7!} +\binom{7}{7}\frac{\lambda^8}{8!}
  \\
  &\quad
  -\lambda \sum_{j=3}^{5} \binom{5}{j}\,\frac{\lambda^j}{j!}
  +\frac{\lambda^2}{2!} \sum_{j=2}^{3} \binom{3}{j+1}\,\frac{\lambda^j}{j!}
  -\frac{\lambda^3}{3!} \sum_{j=1}^{1} \binom{1}{j+2}\,\frac{\lambda^j}{j!}
  \\
  &= \frac{35\lambda^4}{4!} +\frac{35\lambda^5}{5!} +\frac{21\lambda^6}{6!} +\frac{7\lambda^7}{7!} +\frac{\lambda^8}{8!}
  \\
  &\quad
  -\lambda \Bigl(\frac{10\lambda^3}{3!} +\frac{5\lambda^4}{4!} +\frac{\lambda^5}{5!}\Bigr)
  +\frac{\lambda^4}{2!2!} - 0
  \\
  &= \frac{\lambda^4}{0!4!} +\frac{\lambda^5}{2!3!} +\frac{\lambda^6}{4!2!} +\frac{\lambda^7}{6!1!} +\frac{\lambda^8}{8!0!} \,.
\end{split}
\eq
\end{enumerate}


\end{document}